\input amstex.tex
\input amsppt.sty
\input epsf.sty

\define\la{\lambda}
\define\La{\Lambda}
\define\Y{\Bbb Y}
\define\Z{\Bbb Z}
\redefine\S{\text{\bf S}}
\define\M{\text{\bf M}}
\define\T{\text{\bf T}}
\define\SS{\Cal S}
\define\C{\Bbb C}

\define\lam{\vec{\lambda}}
\redefine\ss{\Cal X}
\define\sss{\Cal Y}

\define\tht{\thetag}

\NoBlackBoxes

\topmatter

\title Schur dynamics of the Schur processes
\endtitle

\author Alexei Borodin
\endauthor

%\date Preliminary version, January 18, 2010
%\enddate

\abstract We construct discrete time Markov chains that preserve the
class of Schur processes on partitions and signatures.

One application is a simple exact sampling algorithm for
$q^{volume}$-distributed skew plane partitions with an arbitrary
back wall. Another application is a construction of  Markov chains
on infinite Gelfand-Tsetlin schemes that represent deterministic
flows on the space of extreme characters of the infinite-dimensional
unitary group.
\endabstract

\toc \widestnumber\head{10} \head {} Introduction\endhead

\head 1. Nonnegative specializations of the Schur functions
\endhead

\head 2. The Schur processes
\endhead

\head 3. Example 1. Measures $q^{volume}$ on skew plane partitions
\endhead

\head 4. Example 2. Path measures for extreme characters of
$U(\infty)$
\endhead

\head 5. Markov chains on the Schur processes
\endhead

\head 6. Markov chains on the two-sided Schur processes
\endhead

\head 7. Exact sampling algorithms
\endhead

\head 8. A general construction of multivariate Markov chains
\endhead

\head 9. Application to the Schur processes
\endhead

\head 10. Application to the two-sided Schur processes\endhead

\head{} References \endhead

\endtoc

\endtopmatter

\document

\head Introduction \endhead

The {\it Schur processes\/} were introduced in \cite{OR1} as a class
of measures on sequences of partitions in order to study large
random plane partitions with weights proportional to $q^{volume}$,
$0<q<1$. The concept generalized that of the {\it Schur measures\/}
introduced earlier in \cite{Ok}. The asymptotic techniques of
\cite{OR1} were developed further in \cite{OR2} to study the
asymptotics of large skew plane partitions, see also \cite{BMRT}.

The range of applications of the Schur measures and Schur processes
expanded quickly; apart from random plane partitions they have been
applied to harmonic analysis on the infinite symmetric group
\cite{Ok}, Szeg\"o-type formulas for Toeplitz determinants
\cite{BO}, relative Gromov-Witten theory of $\C^*$ \cite{OP}, random
domino tilings of the Aztec diamond \cite{J2}, discrete and
continuous polynuclear growth processes in 1+1 dimensions \cite{PS},
\cite{J1}, topological string theory \cite{ORV}, and so forth.

The goal of this paper is to define discrete time Markov chains that
map Schur processes to themselves, possibly modifying the
parameters. We also define Markov chains on the {\it two-sided Schur
processes\/} introduced below; the principal difference of those
from the Schur processes is that they live on sequences of {\it
signatures\/} that, unlike partitions, may have negative parts.

The dynamics we construct is also `Schur like'; for example, an
evolution of a partition or a signature that represents a fixed
slice of the (possibly two-sided) Schur process is also a (possibly
two-sided) Schur process.

We present two applications of the construction.

First, we give an exact sampling algorithm for measures of type
$q^{volume}$ on skew plane partitions. Other sampling algorithms for
such measures are known, see \cite{BFP} and references therein.
However, it seems that the algorithm we suggest is simpler; for skew
plane partitions with support fitting in $A\times B$ box, the
algorithm consists in sampling no more that $AB(B+1)/2$ independent
one-dimensional geometric distributions. A short `code' for the
algorithm can be found in Section 7. Exact sampling algorithms for
{\it boxed\/} plane partitions based on similar ideas were
constructed in \cite{BG}, \cite{BGR}.

The second application is a construction of Markov chains on
infinite Gelfand-Tsetlin schemes that preserve the class of Fourier
transforms of the extreme characters of the infinite-dimensional
unitary group, see Section 4 for details. For similar developments
on the infinite-dimensional orthogonal group see \cite{BK}.

A special case of the Markov dynamics that we construct has been
studied in detail in \cite{BF}. One of the goals of this paper is to
provide a more general setup (a broad class of initial conditions
and a multi-parameter family of update rules) for large time
asymptotic analysis of the dynamics.

The construction below is based on a formalism developed in
\cite{BF}, which in its turn was based on an idea from \cite{DF}.
However, our exposition is self-contained.

\subhead Acknowledgements\endsubhead This work was supported in part
by the NSF grant DMS-0707163.

\head 1. Nonnegative specializations of the Schur functions
\endhead

In what follows we use the notation of \cite{M}.

Let $\La$ be the algebra of symmetric functions. A {\it
specialization\/} $\rho$ of $\La$ is an algebra homomorphism of
$\La$ to $\C$; we denote the application of $\rho$ to $f\in\La$ as
$f(\rho)$. The {\it trivial\/} specialization $\varnothing$ takes
value 1 at the constant function $1\in\La$ and takes value $0$ at
any homogeneous $f\in\La$ of degree $\ge 1$.

For two specializations $\rho_1$ and $\rho_2$ we define their union
$\rho=(\rho_1,\rho_2)$ as the specialization defined on Newton power
sums via
$$
p_n(\rho_1,\rho_2)=p_n(\rho_1)+p_n(\rho_2), \qquad n\ge 1.
$$

\example{Definition 1} We say that a specialization $\rho$ of $\La$
is {\it nonnegative\/} if it takes nonnegative values on the Schur
functions: $s_\la(\rho)\ge 0$ for any partition $\la$.
\endexample

The classification of all nonnegative specializations is a classical
result proved independently by Aissen, Edrei, Schoenberg, and
Whitney \cite{AESW} (see also \cite{E}) and Thoma \cite{T}. It says
that a specialization $\rho$ is nonnegative if and only if the
generating function of the images of complete homogeneous functions
has the form
$$
H(\rho;u):=\sum_{n=0}^\infty h_n(\rho) u^n =e^{\gamma u}
\prod_{i\ge1}\frac{1+\beta_i u}{1-\alpha_i u} \tag 1
$$
for certain nonnegative $\{\alpha_i\}$, $\{\beta_i\}$, and $\gamma$
such that $\sum_i(\alpha_i+\beta_i)<\infty$.

It turns out that nonnegativity of  $s_\la(\rho)$ for all $\lambda$
is equivalent to nonnegativity of the images of the skew Schur
functions $s_{\la/\mu}(\rho)$ for all $\la$ and $\mu$. Hence, via
the Jacobi-Trudi formula
$$
s_{\la/\mu}=\det[h_{\la_i-i-\mu_j+j}]_{i,j=1}^r,\qquad r\ge
\max\{\ell(\lambda),\ell(\mu)\},
$$
the classification of nonnegative specializations is equivalent to
that of totally nonnegative triangular Toeplitz matrices with
diagonal entries equal to 1. An excellent exposition of deep
relations of this classification result to representation theory of
the infinite symmetric group can be found in Kerov's book \cite{K}.

For a single $\alpha$ or a single $\beta$ specialization, the values
of skew Schur functions are easy to compute:
$$
\align &H(\rho;u)=\frac1{1-\alpha u}\quad \text{implies}\quad
s_{\la/\mu}(\rho)=\cases
\alpha^{|\la|-|\mu|},&\la_1\ge\mu_1\ge\la_2\ge\mu_2\ge\dots,\\0,&\text{otherwise};
\endcases
\tag 2\\
&H(\rho;u)={1+\beta u}\quad \text{implies}\quad
s_{\la/\mu}(\rho)=\cases \beta^{|\la|-|\mu|},&\la_j-\mu_j\in\{0,1\}\
\text{for all}\ j\ge 1,\\0,&\text{otherwise}.
\endcases
\tag 3
\endalign
$$

 We say that a nonnegative specialization
$\rho$ of $\La$ is {\it admissible} if the generating function
\tht{1} is holomorphic in a disc $D_r=\{u\in\C\mid |u|<r\}$ with
$r>1$. In other words, $\rho$ is admissible iff $\alpha_i<r^{-1}<1$
for all $i$.

Since $H(\rho_1,\rho_2;u)=H(\rho_1;u)H(\rho_2;u)$, the union of
admissible specializations is admissible (unions of nonnegative specializations
are also nonnegative).

For a nonnegative specialization $\rho$, denote by $\Y(\rho)$ the
set of partitions (or Young diagrams) $\lambda$ such that
$s_\lambda(\rho)>0$. We also call $\Y(\rho)$ the {\it support\/} of
$\rho$. The set of all partitions will be denoted as $\Y$.

Using the combinatorial formula for the Schur functions \cite{M,
Sect. I.5 (5.12)} and the involution $\omega$ \cite{M, Sect. I.2},
it is not hard to show that if, for a nonnegative specialization
$\rho$, in \tht{1} $\gamma=0$ and there are $p<\infty$ nonzero
$\alpha_j$'s and $q<\infty$ nonzero $\beta_j$'s, then $\Y(\rho)$
consists of the Young diagrams that fit into the $\Gamma$-shaped
figure with $p$ rows and $q$ columns. Otherwise it is easy to see
that $\Y(\rho)=\Y$.

In particular, if in \tht{1} all $\beta_j$'s and $\gamma$ vanish,
and there are $p$ nonzero $\alpha_j$'s, then $\Y(\rho)$ consists of
Young diagrams with no more than $p$ rows. Such a specialization
consists in assigning values $\alpha_j$ to $p$ of the symmetric
variables used to define $\La$, and $0$'s to all the other symmetric
variables.

We will also need minors of arbitrary (not necessarily triangular)
doubly-infinite totally nonnegative Toeplitz matrices. The
classification of such matrices was obtained by Edrei in \cite{E},
who proved an earlier conjecture of Schoenberg. The result is as
follows.

A matrix $M=[M_{i-j}]_{i,j=-\infty}^{+\infty}$ is totally
nonnegative if and only if, after a transformation  of the form
$M_n\mapsto cR^nM_n$ with $c>0$, $R\ge 0$, the generating function
of its entries has the form
$$
\multline H(M;u):=\sum_{n=-\infty}^{+\infty} M_n u^n\\
=e^{\gamma^+(u-1)+\gamma^-(u^{-1}-1)}\prod_{i=1}^\infty\left(
\frac{1+\beta_i^+(u-1)}{1-\alpha_i^+(u-1)}\,\frac{1+\beta_i^-(u^{-1}-1)}{1-\alpha_i^-(u^{-1}-1)}
\right)
\endmultline
\tag 4
$$
for certain nonnegative $\{\alpha_j^\pm\}$, $\{\beta_j^\pm\}$, and
$\gamma^\pm$ such that
$\sum(\alpha_i^++\alpha_i^-+\beta_i^++\beta_i^-)<\infty$ and
$\beta_j\le 1$ for all $j$. The parametrization of $M$ by
$(\{\alpha_j^\pm\},\,\{\beta_j^\pm\},\,\gamma^\pm)$ becomes unique
if one adds the condition
$\max_j\{\beta_j^+\}+\max_j\{\beta_j^-\}\le 1$.

The generating function on the left is understood as the Laurent
series of the holomorphic function in a neighborhood of the unit
circle $|u|=1$ that stands on the right. We call the largest annulus
of the form $\{u\in\C\mid 0\le r_1<|u|<r_2\}$ where $H(M;u)$ is
holomorphic (the unit circle must be inside the annulus) the {\it
analyticity annulus\/} of $H(M;u)$.

\example{Definition 2} We say that a totally nonnegative Toeplitz
matrix $M$ is {\it admissible\/} if the generating function of its
entries is given by \tht{4} (i.e., no multiplication by $cR^n$ is
involved).
\endexample

Note that since multiplying Toeplitz matrices corresponds to
multiplying the generating functions \tht{4}, the product of two
admissible matrices is admissible.

It will be convenient for us to use a similar notation for the
minors of general Toeplitz matrices as in the triangular case (Jacobi-Trudi
formula).

Define {\it signatures\/} of length $n$ as $n$-tuples
$\la=(\la_1\ge\la_2\ge\dots\ge \la_n)$ of nonincreasing integers. We
will also write $\ell(\la)=n$ and $|\la|=\la_1+\la_2+\dots+\la_n$.
By convention, there is a unique signature $\varnothing$ of length
$0$ with $|\varnothing|=0$.

For any two signatures $\la$ and $\mu$ of length $n$ and an
admissible $M$ we set
$$
s_{\la/\mu}(M)=\det\left[M_{\la_i-i-\mu_j+j}\right]_{i,j=1}^n.
$$
For totally nonnegative $M$ with only one $\alpha^{\pm}$ or
$\beta^\pm$ parameter nonzero (and all other parameters being zero),
one obtains formulas analogous to \tht{2}, \tht{3}:
$$
H(\rho;u)=\frac1{1-\alpha (u^{\pm 1}-1)}\quad \text{implies}\quad
s_{\la/\mu}(\rho)=\dfrac
1{(1+\alpha)^{n}}\left(\dfrac{\alpha}{1+\alpha}\right)^{\pm|\la|\mp|\mu|}
\tag 5
$$
if $\pm\la_j\mp\mu_j\ge 0$ for all $1\le j\le n$, and 0 otherwise;
$$
H(\rho;u)={1+\beta (u^{\pm 1}-1)}\quad  \text{implies} \quad
s_{\la/\mu}(\rho)=(1-\beta)^n\left(\frac\beta{1-\beta}\right)^{\pm|\la|\mp|\mu|}\tag
6
$$
if $\pm\la_j\mp\mu_j\in\{0,1\}$ for all $1\le j\le n$, and 0
otherwise.

Also, mimicking the property of the Schur functions, for a constant
$c\in\C$, a signature $\nu$ of length $n+1$, and a signature $\la$
of length $n$, we set
$$
s_{\la/\mu}(c):=\cases c^{|\la|-|\mu|}, & \la_{n+1}\le \mu_{n}\le
\la_{n}\le \dots\le \la_2 \le \mu_1\le \la_1,\\
0,&\text{otherwise},
\endcases
\tag 7
$$
with the convention that $0^0=1$.

\head 2. The Schur processes
\endhead

Pick a natural number $N$ and admissible specializations
$\rho_0^+,\dots,\rho_{N-1}^+$, $\rho_1^-,\dots,\rho_N^-$ of $\La$.
For any sequences $\la=(\la^{(1)},\dots,\la^{(N)})$ and
$\mu=(\mu^{(1)},\dots,\mu^{(N-1)})$ of partitions satisfying
$$
\varnothing\subset \la^{(1)}\supset \mu^{(1)}\subset
\la^{(2)}\supset \mu^{(2)} \subset \dots \supset \mu^{(N-1)}\subset
\la^{(N)}\supset \varnothing \tag 8
$$
define their weight as
$$
\Cal
W(\la,\mu):=s_{\la^{(1)}}(\rho_0^+)\,s_{\la^{(1)}/\mu^{(1)}}(\rho_1^-)
s_{\la^{(2)}/\mu^{(1)}}(\rho_1^+)\,\cdots
s_{\la^{(N)}/\mu^{(N-1)}}(\rho_{N-1}^+)\, s_{\la^{(N)}}(\rho_N^-).
\tag 9
$$

There is one Schur function factor for any two neighboring
partitions in \tht{8}.

The fact that all the specializations are nonnegative implies that
all the weights are nonnegative. The admissibility of $\rho$'s
implies that
$$
Z(\rho_0^+,\dots,\rho_{N-1}^+;\rho_1^-,\dots,\rho_N^-):=\sum_{\la,\mu}
\Cal W(\la,\mu)=\prod_{0\le i<j\le N}H(\rho_i^+;\rho_j^-)<\infty,
\tag 10
$$
where $H(\rho_1;\rho_2)=\sum_{\la\in\Y}s_\la(\rho_1)s_\la(\rho_2)=
\exp\left(\sum_{n\ge 1}p_n(\rho_1)p_n(\rho_2)/n\right)$, and $p_n$'s
are the Newton power sums. Indeed, this follows from the repeated
use of identities, cf. \cite{M, I(5.9) and Ex.~I.5.26(1)},
$$
\gather \sum_{\kappa\in\Y}
s_{\kappa/\nu}(\rho_1)s_{\kappa/\hat\nu}(\rho_2)=H(\rho_1;\rho_2)\sum_{\tau\in\Y}
s_{\nu/\tau}(\rho_2)s_{\hat\nu/\tau}(\rho_1),\tag 11\\
\sum_{\nu\in\Y}s_{\kappa/\nu}(\rho_1)s_{\nu/\tau}(\rho_2)=s_{\kappa/\tau}(\rho_1,\rho_2),
\tag 12
\endgather
$$
and from the fact that for an admissible specialization $\rho$ with $H(\rho;u)$
holomorphic in a disc of radius $r$, we have
$p_n(\rho)=O(r^{-n})$.

The same argument shows that the partition function \tht{10} is
finite under the weaker assumption of finiteness of all
$H(\rho_i^+;\rho_j^-)$ for $0\le i<j\le N$.

\example{Definition 3} The {\it Schur process\/}
$\S(\rho_0^+,\dots,\rho_{N-1}^+;\rho_1^-,\dots,\rho_N^-)$ is the
probability distribution on sequences $(\la,\mu)$ as in \tht{8} with
$$
\S(\rho_0^+,\dots,\rho_{N-1}^+;\rho_1^-,\dots,\rho_N^-)(\la,\mu)=\frac{\Cal
W(\la,\mu)}
{Z(\rho_0^+,\dots,\rho_{N-1}^+;\rho_1^-,\dots,\rho_N^-)}\,.
$$
The Schur process with $N=1$ is called the {\it Schur measure\/}.
\endexample

Using \tht{11}-\tht{12} it is not difficult to show that a
projection of the Schur process to any subsequences of $(\la,\mu)$
is a also a Schur process. In particular, the projection of
$\S(\rho_0^+,\dots,\rho_{N-1}^+;\rho_1^-,\dots,\rho_N^-)$ to
$\lambda^{(j)}$ is the Schur measure
$\S(\rho^+_{[0,j-1]};\rho^-_{[j,N]})$, and its projection to
$\mu^{(k)}$ is a slightly different Schur measure
$\S(\rho^+_{[0,k-1]};\rho^-_{[k+1,N]})$. Here we used the notation
$\rho^\pm_{[a,b]}$ to denote the union of specializations
$\rho^\pm_m$, $m=a,\dots,b$.

We now aim at defining a Schur like process for signatures.

Pick a natural number $N$, real numbers  $a_1,\dots,a_N>0$,
nonnegative integers $c(1),\dots,c(N)$, and $c(1)+\dots+c(N)$
admissible Toeplitz matrices
$$
\M=\{M^{(k,l)}\mid 1\le k\le N,\, 1\le l\le c(k)\}.
$$
If all $c(k)$ are zero then $\M$ is empty.

We will also need a totally nonnegative matrix of size $\Z \times
N$, denote it as
$$
\Psi={[\Psi_{ij}]}_{i\in\Z,-1\ge j\ge -N}.
$$

For any sequences
$$
\lam^{(1)}=(\la^{(1,0)},\la^{(1,1)},\dots,\la^{(1,c(1))}),\,\dots,\,
\lam^{(N)}=(\la^{(N,0)},\la^{(N,1)},\dots,\la^{(N,c(N))}) \tag 13
$$
of signatures of lengths $\ell(\la^{(k,*)})=k$, define their
(nonnegative) weight as
$$
\multline \Cal W(\lam^{(1)},\dots,\lam^{(N)}):=
\det\left[\Psi_{\lambda_i^{(N,c(N))}-i,-j}\right]_{i,j=1}^N\\ \times
\prod_{k=1}^N \left(s_{\lambda^{(k,0)}/\lambda^{(k-1,c(k-1))}}(a_k)
\prod_{l=1}^{c(k)}
s_{\lambda^{(k,l)}/\lambda^{(k,l-1)}}\left(M^{(k,l)}\right)\right)
\endmultline
\tag 14
$$
with $\lambda^{(0,c(0))}=\varnothing$.

We assume that the generating functions
$$
\Psi_j(u):=\sum_{n=-\infty}^{+\infty} \Psi_{n,-j}\,u^{n+j} \tag 15
$$
are holomorphic in an open set containing the unit circle. As we
will see in Section 10, if for any $j\le N$, $a_j$ lies in the
common analyticity annulus for $\{H(M^{(k,l)};u^{-1})\}_{k\ge j}$,
$\{\Psi_i(u)\}_{i=1}^N$, then the partition function of weights \tht{14}
is finite and it has the form
$$
\multline
Z(a_1,\dots,a_N;\M;\Psi):=\sum_{\lam^{(1)},\dots,\lam^{(N)}}\Cal
W(\lam^{(1)},\dots,\lam^{(N)})\\ =
\frac{\det\left[a_i^{-j}\Psi_j(a_i)\right]_{i,j=1}^N}{\det\left[a_i^{-j}\right]_{i,j=1}^N}\,
\prod\limits_{1\le
j \le k\le
N}\prod\limits_{l=1}^{c(k)}H\left(M^{(k,l)};a_j^{-1}\right)\,.
\endmultline
\tag 16
$$

In the important special case when the matrix $\Psi$ is actually
Toeplitz, $\Psi_{i,-j}=\psi_{i+j}$, \tht{16} simplifies:
$$
Z(a_1,\dots,a_N;\M;\Psi)=\prod_{i=1}^N \psi(a_i)\,\prod\limits_{1\le
j \le k\le
N}\prod\limits_{l=1}^{c(k)}H\left(M^{(k,l)};a_j^{-1}\right), \tag 17
$$
where $\psi(u)=\sum_{n\in\Z}\psi_n u^n$.

\example{Definition 4} The {\it two-sided Schur process\/}
$\T(a_1,\dots,a_N;\M;\Psi)$ is the probability distribution on
sequences $(\la^{(1)},\dots,\la^{(N)})$ as in \tht{13} with
$$
\T(a_1,\dots,a_N;\M;\Psi)(\lam^{(1)},\dots,\lam^{(N)})=\frac{\Cal
W(\lam^{(1)},\dots,\lam^{(N)})} {Z(a_1,\dots,a_N;\M;\Psi)}\,.
$$
\endexample

\example{Remark 5} If in the Schur process of Definition 3 each of
the specializations $\rho_j^+$ is a one-variable specialization with
$H(\rho_j^+;u)=(1-a_{j+1}u)^{-1}$, $j=0,\dots,N-1$, then the Schur
process can be viewed as a special case of the two-sided Schur
process with $c(1)=\dots=c(N-1)=1$, $c(N)=0$, and identification
$$
\gathered \lambda^{(j)}=\lambda^{(j,0)}, \quad j=1,\dots,N,\qquad
\mu^{(j)}=\lambda^{(j,1)},\quad j=1,\dots,N-1,\\
H(\rho_k^-;u)=H(M^{(k,1)};u^{-1}),\quad k=1,\dots,N-1; \qquad
H(\rho_N^-;u)=\psi(u).
\endgathered
$$
The corresponding two-sided Schur process lives on signatures with
nonnegative parts that can also be viewed as partitions.

Observe that under this identification the formulas \tht{10} and
\tht{17} coincide.
\endexample

\head 3. Example 1. Measures $q^{volume}$ on skew plane partitions
\endhead

Fix two natural numbers $A$ and $B$. For a Young diagram $\pi\subset
B^A$, set $\bar\pi=B^A/\pi$.

A (skew) plane partition $\Pi$ with support $\bar\pi$ is a filling
of all boxes of $\bar\pi$ by nonnegative integers $\Pi_{i,j}$ (we
assume that $\Pi_{i,j}$ is located in the $i$th row and $j$th column
of $B^A$) such that $ \Pi_{i,j}\ge \Pi_{i,j+1}$ and $\Pi_{i,j}\ge
\Pi_{i+1,j}$ for all values of $i,j$.

The volume of the plane partition $\Pi$ is defined as
$$
vol(\Pi)=\sum_{i,j}\Pi_{i,j}.
$$

The goal of the section is to explain that the measure on plane
partitions with given support $\bar\pi$  and weights proportional to
$q^{vol(\,\cdot\,)}$, $0<q<1$, is a Schur process. This fact has
been observed and used in \cite{OR1}, \cite{OR2}, \cite{BMRT}.

The Schur process will be such that for any two neighboring
specializations $\rho_k^-,\rho_k^+$ at least one is trivial. This
implies that each $\mu^{(j)}$ coincides either with $\lambda^{(j)}$
or with $\lambda^{(j+1)}$. Thus, we can restrict our attention to
$\lambda^{(j)}$'s only.

 For a plane partition $\Pi$, we set ($1\le k\le A+B+1$)
$$
\lambda^{(k)}(\Pi)=\bigl\{\Pi_{i,i+k-A-1}\mid
(i,i+k-A-1)\in\bar\pi\bigr\}.
$$
Note that $\lambda^{(1)}=\lambda^{(A+B+1)}=\varnothing$.

We need one more piece of notation. Define
$$
\Cal L(\pi)=\{A+\pi_i-i+1\mid i=1,\dots,A\}.
$$
This is an $A$-point subset in $\{1,2,\dots,A+B\}$, and all such
subsets are in bijection with the partitions $\pi$ contained in the
box $B^A$. The elements of $\Cal L(\pi)$ mark the ``up-steps'' in
the boundary of $\pi$ (=back wall of $\Pi$).

%\vskip 0.5cm\hskip 1cm
 %{\pdfximage height 40 mm
%{plane_partition.pdf}\pdfrefximage \pdflastximage}

%\vskip -4.5cm\hskip 5.8cm {\pdfximage width 55 mm
%{plane_partition.jpg}\pdfrefximage \pdflastximage}

$$ \epsffile{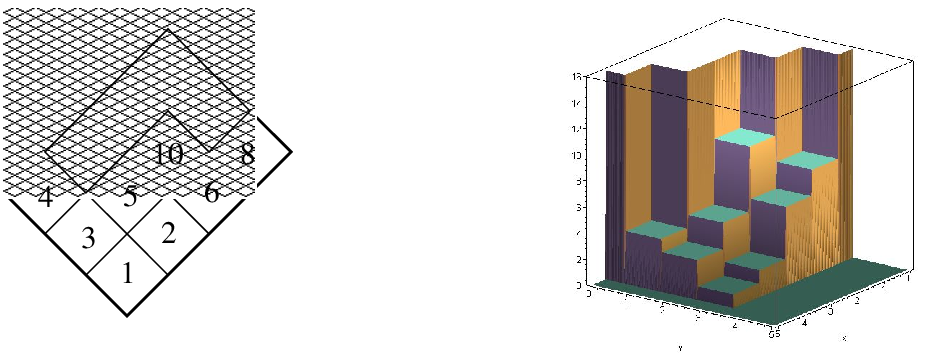} $$

The figure above shows a plane partition $\Pi$ and its plot with
$$
\gathered A=4,\ B=3,\ \pi=(2,1,1,0),\\
\lambda^{(2)}=(4),\ \lambda^{(3)}=3,\ \lambda^{(4)}=(5,1),\
\lambda^{(5)}=(10,2),\ \lambda^{(6)}=(6),\ \lambda^{(7)}=(8),\\
vol(\Pi)=\sum_{i=2}^{A+B}|\la^{(i)}|=39,\quad \Cal
L(\pi)=\{1,3,4,6\}.
\endgathered
$$

\proclaim{Proposition 6} Let $\pi$ be a partition contained in the
box $B^A$. The measure on the plane partitions $\Pi$ with support
$\bar\pi$ and weights proportional to $q^{vol(\Pi)}$, is the Schur
process with $N=A+B+1$ and nonnegative specializations
$\{\rho_i^+\}$, $\{\rho_j^-\}$ defined by
$$
\gather H(\rho_0^+;u)=H(\rho_N^-;u)=1,\\
H(\rho_j^+;u)=\cases \dfrac 1{1-q^{-j}u},&j\in\Cal L(\pi),
\\ 1,&j\notin\Cal L(\pi);
\endcases
\qquad H(\rho_j^-;u)=\cases 1,&j\in\Cal L(\pi),\\\dfrac
1{1-q^{j}u},&j\notin\Cal L(\pi).
\endcases
\endgather
$$
\endproclaim

Note that not all specializations are admissible, but the weaker
assumption of finiteness of $H(\rho_i^+;\rho_j^-)$ for $0\le i<j\le
N$ guarantees that the partition function is finite.

\demo{Proof} Observe that the set of all plane partitions supported
by $\bar\pi$, as well as the support of the Schur process from the
statement of the proposition, consists of sequences
$(\lambda^{(1)},\lambda^{(2)},\dots,\lambda^{(N)})$ with
$$
\gathered \lambda^{(1)}=\lambda^{(N)}=\varnothing, \\
\lambda^{(j)}\prec\lambda^{(j+1)} \text{  if  } j\in \Cal
L(\lambda),\qquad \lambda^{(j)}\succ \lambda^{(j+1)} \text{  if  }
j\notin \Cal L(\lambda),
\endgathered
$$
where we write $\mu\prec\nu$ or $\nu\succ\mu$  iff
$\nu_1\ge\mu_1\ge\nu_2\ge\mu_2\ge \dots$\,.

On the other hand, \tht{2} implies that the weight of
$(\lambda^{(1)},\lambda^{(2)},\dots,\lambda^{(N)})$ with respect to
the Schur process from the hypothesis is equal to $q$ raised to the power
$$
\sum_{j=2}^{A+B}|\lambda^{(j)}|\Bigl( -(j-1){\bold 1}_{j-1\in \Cal
L(\pi)}-(j-1){\bold 1}_{j-1\notin \Cal L(\pi)} +j{\bold 1}_{j\in
\Cal L(\pi)}+j{\bold 1}_{j\notin \Cal L(\pi)} \Bigr),
$$
where the four terms are the contributions of
$\rho_{j-1}^+,\rho_{j-1}^-,\rho_j^+,\rho_j^-$, respectively.

Clearly, the sum is equal to
$\sum_{j=2}^{A+B}|\lambda^{(j)}|=vol(\Pi)$. \qed
\enddemo

\example{Remark 7} A similar statement holds for any measure on
plane partitions with weights proportional to $\prod q_j^{|\la_j|}$
with possibly different positive parameters $q_j$, as long as the
partition function is finite. The proof is very similar.
\endexample

\head 4. Example 2. Path measures for extreme characters of
$U(\infty)$
\endhead

Let $U(N)$ denote the group of $N\times N$ unitary matrices. It is a
classical result that the irreducible representations of $U(N)$ can
be paramterized by signatures
$\lambda=(\lambda_1\geq\ldots\geq\lambda_N)$ of length $N$ also
called {\it highest weights\/}. Thus, there is a natural bijection
$\lambda\leftrightarrow\chi^{\lambda}$ between signatures of length
$N$ and the conventional irreducible characters (=traces of
irreducible representations) of $U(N)$.

For each $N$, embed $U(N)$ in $U(N+1)$ as the subgroup fixing the
$(N+1)$st basis vector. Equivalently, each $U\in U(N)$ can be
thought of as an $(N+1)\times(N+1)$ matrix by setting
$U_{i,N+1}=U_{N+1,j}=0$ for $1\leq i,j\leq N$ and $U_{N+1,N+1}=1$.
The union $\bigcup_{N=1}^{\infty} U(N)$ is denoted $U(\infty)$ and
called the {\it infinite-dimensional unitary group\/}.

A {\it character\/} of $U(\infty)$ is a positive definite function
$\chi:U(\infty)\rightarrow\C$ which is constant on conjugacy classes
and normalized by $\chi(\bold{1})=1$. We further assume that $\chi$
is continuous on each $U(N)\subset U(\infty)$. The set of all
characters of $U(\infty)$ is convex, and the extreme points of this
set are called {\it extreme\/} characters.

Remarkably, the extreme characters of $U(\infty)$ are in one-to-one
correspondence with admissible Toeplitz matrices $M$ from Definition
2, see \cite{Vo}, \cite{VK}, \cite{OO}. The values of the character
$\chi^M$ corresponding to $M$ are given by
$$
\chi^M(U)=\prod_{u\in \operatorname{Spectrum}(U)} H(M;u),
$$
where $H(M;u)$ is given in \tht{4}.

Let $\Bbb{GT}_N$ be the set of all signatures of length $N$; set
$\Bbb{GT}=\bigcup_N\Bbb{GT}_N$. Turn $\Bbb{GT}$ into a graph by
drawing an edge between signatures $\lambda\in\Bbb{GT}_N$ and
$\mu\in\Bbb{GT}_{N+1}$ if $\lambda$ and $\mu$ satisfy the branching
relation $\lambda\prec\mu$, where $\lambda\prec\mu$ means that
$\mu_1\leq\lambda_1\leq\mu_2\leq\lambda_2\leq\ldots\leq\lambda_N\leq\mu_{N+1}$.
$\Bbb{GT}$ is known as the {\it Gelfand-Tsetlin graph\/}.

A {\it path} in $\Bbb{GT}$, or an {\it infinite Gelfand-Tsetlin
scheme\/}, is an infinite sequence $t=(t_1,t_2,\ldots)$ such that
$t_i\in\Bbb{GT}_i$ and $t_i\prec t_{i+1}$. Let $\Cal T$ be the set
of all such paths.

One can also look at finite paths, or finite Gelfand-Tsetlin
schemes, which are sequences $\tau=(\tau_1,\tau_2,\ldots,\tau_N)$
such that $\tau_i\in\Bbb{GT}_i$ and
$\tau_1\prec\tau_2\prec\ldots\prec\tau_N$. Denote the set of all
paths of length $N$ by $\Cal T_N$.

%\vskip 0.5cm \hskip 2cm
% {\pdfximage height 40 mm
%{gt_scheme.pdf}\pdfrefximage \pdflastximage}
%\vskip 0.2cm

%\vskip -4.8cm\hskip 5.8cm {\pdfximage width 55 mm
%{gt_scheme.jpg}\pdfrefximage \pdflastximage}

$$ \epsffile{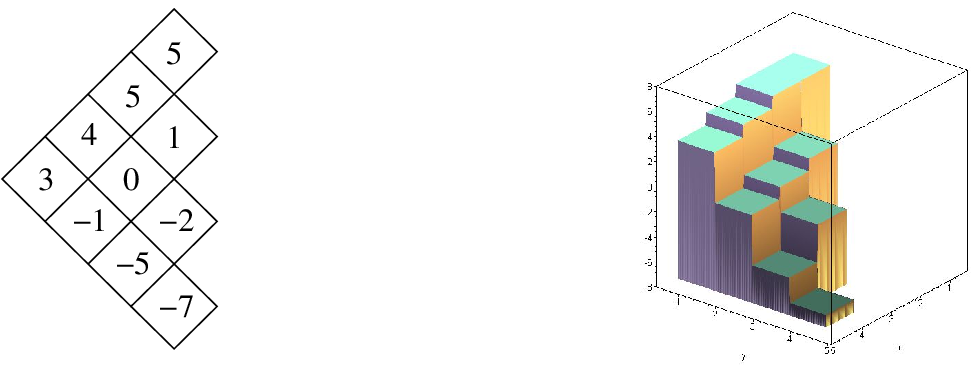} $$

The figure above depicts a Gelfand-Tsetlin scheme $\tau\in\Cal T_4$
and its plot with
$$
\tau_1=(3),\quad \tau_2=(4,-1),\quad \tau_3=(5,0,-5),\quad
\tau_4=(5,1,-2,-7).
$$

Each character $\chi$ of $U(\infty)$ defines a probability measure
$P_N^{\chi}$ on $\Bbb{GT}_N$: Restricting the character to $U(N)$,
we have
$$
\chi\Bigl|_{U(N)} = \sum_{\lambda\in\Bbb{GT}_N}
P_N^\chi(\lambda)\frac{{\chi}^{\lambda}}{{\chi}^{\lambda}(\bold{1}_N)}\,.
$$

For each finite path $\tau\in\Cal T_N$, let $C_{\tau}\subset\Cal T$
be the set
$$
C_{\tau}=\{t\in\Cal T:(t_1,t_2,\ldots,t_N)=\tau\}.
$$

A character $\chi$ of $U(\infty)$ also defines a probability measure
$P^{\chi}$ on $\Cal T$ (with a suitably defined Borel structure),
which can be uniquely specified by setting
$$P^{\chi}(C_{\tau})=\frac{P_N^{\chi}(\lambda)}{\chi^\la(\bold{1}_N)},$$
where $\tau$ is an arbitrary finite path ending at $\lambda$, see
\cite{Ol, Section 10} for details. Note that we assign the same
weight to all finite paths with the same end.

We use the same formula to define a probability
measure $P^\chi_{[1,N]}$ on $\Cal T_N$, which is just the projection
of $P^\chi$ from $\Cal T$ to $\Cal T_N$.

\proclaim{Proposition 8} For any admissible Toeplitz matrix $M$ as
in Definition 2, the measure $P^{\chi^M}_{[1,N]}$ on $\Cal T_N$
coincides with the two-sided Schur process of Definition 4 with
$$
a_1=\dots=a_N=1,\quad c(1)=\dots=c(N)=0, \quad \Psi=M,
$$
and with sequences $(\la^{(1,0)},\dots,\la^{(N,0)})$ viewed as
elements of $\Cal T_N$.
\endproclaim
\demo{Proof} Directly follows from \tht{7} and Lemma 6.5 of
\cite{Ol}.\qed
\enddemo

\head 5. Markov chains on the Schur processes
\endhead

Let us introduce some notation.

For two nonnegative specializations $\rho_1,\rho_2$ of $\La$ such
that $H(\rho_1;\rho_2)<\infty$, and $\la,\mu\in\Y$, set
$$
P_{\rho_1,\rho_2}(\la,\mu \uparrow\uparrow \nu)=const\cdot
s_{\nu/\la}(\rho_1)s_{\nu/\mu}(\rho_2),\qquad \nu\in\Y,
$$
where we assume that
$$
\{\nu\in\Y\mid
s_{\nu/\la}(\rho_1)s_{\nu/\mu}(\rho_2)>0\}\ne\varnothing, \tag 18
$$
and the constant prefactor is chosen so that we obtain a probability
measure in $\nu$:
$$
\sum_{\nu\in\Y} P_{\rho_1,\rho_2}(\la,\mu \uparrow\uparrow \nu)=1.
$$
Given \tht{18}, the existence of such constant follows from \tht{11}.

Similarly, dropping the assumption $H(\rho_1;\rho_2)<\infty$, we
define
$$
\gather
P_{\rho_1,\rho_2}(\la,\mu \downarrow\uparrow \nu)=const\cdot
s_{\la/\nu}(\rho_1)s_{\nu/\mu}(\rho_2),\\
P_{\rho_1,\rho_2}(\la,\mu \uparrow\downarrow \nu)=const\cdot
s_{\nu/\la}(\rho_1)s_{\mu/\nu}(\rho_2),\\
P_{\rho_1,\rho_2}(\la,\mu \downarrow\downarrow \nu)=const\cdot
s_{\la/\nu}(\rho_1)s_{\mu/\nu}(\rho_2),
\endgather
$$
where in all three cases we assume that the set of $\nu$ giving
nonzero values on the right-hand side is nonempty (it is
finite in all three cases), and we choose constants so that we obtain probability
distributions in $\nu\in\Y$.

If both $\rho_1$ and $\rho_2$ are single-$\alpha$ or single-$\beta$
specializations, relations \tht{2}, \tht{3} show that all four
distributions $P_{\rho_1,\rho_2}$ are products of geometric
distributions conditioned to stay in segments and Bernoulli
measures.

\example{Example 8} Assume that $H(\rho_1;u)=(1-au)^{-1}$,
$H(\rho_2;u)=(1-bu)^{-1}$. Denote by $G_{m,n}^\xi$, $m\le n$, the
probability distribution on the set $\{m,m+1,\dots,n\}$ given by
$$
G_{m,n}^\xi(\{k\})=\frac{\xi^k}{\sum_{j=m}^n
\xi^j}=\frac{1-\xi^{n-m+1}}{\xi^m(1-\xi)}\cdot\xi^{k},\qquad m\le
k\le n.
$$
Then
$$
\align &P_{\rho_1,\rho_2}(\la,\mu \uparrow\uparrow
\nu)=G_{\max\{\la_1,\mu_1\},+\infty}^{ab}(\nu_1)\prod_{j\ge 2}
G_{\max\{\la_j,\mu_j\},\min\{\la_{j-1},\mu_{j-1}\}}^{ab}(\nu_j),\\
& P_{\rho_1,\rho_2}(\la,\mu \downarrow\uparrow
\nu)=G^{b/a}_{\max\{\la_{2},\mu_1\},\la_1}(\nu_1)\prod_{j\ge 2}
G^{b/a}_{\max\{\la_{j+1},\mu_{j}\},\min\{\la_{j},\mu_{j-1}\}}(\nu_j),
\endalign
$$
where in the first case we need to additionally assume that $ab<1$
(equivalently, $H(\rho_1;\rho_2)<\infty$).

Further, assume that $H(\rho_3;u)=(1+cu)$. Denote by $B_{m,n}^p$,
$n\in\{m,m+1\}$, the probability distribution on $\{m,m+1\}$ given
by
$$
B_{m,m}^p(\{k\})=\cases 1,&k=m,\\ 0,& k=m+1,\endcases
\qquad B_{m,m+1}^p(\{k\})=\cases \dfrac1{1+c},&k=m,\\\dfrac{c}{1+c},&
k=m+1.\endcases
$$

Then
$$
\align &P_{\rho_3,\rho_2}(\la,\mu \uparrow\uparrow
\nu)=B_{\max\{\la_1,\mu_1\},\la_1+1}^{bc}(\nu_1)\prod_{j\ge 2}
B_{\max\{\la_j,\mu_j\},\min\{\la_{j-1},\la_{j}+1,\mu_{j-1}\}}^{bc}(\nu_j),\\
& P_{\rho_3,\rho_2}(\la,\mu \downarrow\uparrow
\nu)=B^{b/c}_{\max\{\la_{1}-1,\la_2,\mu_1\},\la_1}(\nu_1)\prod_{j\ge
2}
B^{b/c}_{\max\{\la_{j}-1,\la_{j+1},\mu_{j}\},\min\{\la_{j},\mu_{j-1}\}}(\nu_j).
\endalign
$$
\endexample

Let $(\rho_0^+,\dots,\rho_{N-1}^+;\rho_1^-,\dots,\rho_N^-)$ be
nonnegative specializations of $\La$ defining a Schur process as in
Definition 3. Let $\pi$ be another nonnegative specialization of
$\La$ such that $H(\pi,\rho^+_j)<\infty$ for all $0\le j<N$.

Let $\ss$ be the set of pairs of sequences $(\la,\mu)$ as in \tht{8}
with
$$
s_{\la^{(1)}}(\rho_0^+)\,s_{\la^{(1)}/\mu^{(1)}}(\rho_1^-)
s_{\la^{(2)}/\mu^{(1)}}(\rho_1^+)\,\cdots
s_{\la^{(N)}/\mu^{(N-1)}}(\rho_{N-1}^+)>0
$$
The product above is the same as in \tht{9} without the last factor.
Thus, the support of
$\S(\rho_0^+,\dots,\rho_{N-1}^+;\rho_1^-,\dots,\rho_N^-)$ is
contained in $\ss$.

Define a matrix $\frak P^\uparrow_\pi$ with rows and columns
parameterized by elements of $\ss$ via
$$
\multline \frak P^\uparrow_\pi((\la,\mu),(\tilde\la,\tilde\mu))=
P_{\rho_0^+,\pi}\bigl(\varnothing,\la^{(1)}\uparrow\uparrow
\tilde\la^{(1)}\bigr)\\ \times \prod_{j=1}^{N-1}
P_{\rho_j^-,\pi}\bigl(\tilde\la^{(j)},\mu^{(j)}\downarrow\uparrow
\tilde\mu^{(j)}\bigr)P_{\rho_j^+,\pi}\bigl(\tilde\mu^{(j)},\la^{(j+1)}\uparrow\uparrow
\tilde\la^{(j+1)}\bigr).
\endmultline
\tag 19
$$

In other words, starting from $(\la,\mu)$, one first finds
$\tilde\la^{(1)}$ using $\la^{(1)}$, then $\tilde\mu^{(1)}$ using
$\tilde\la^{(1)}$ and $\mu^{(1)}$, then $\tilde\la^{(2)}$ using
$\tilde\mu^{(1)}$ and $\la^{(2)}$, and so on. One could say that we
perform {\it sequential update\/}.

Note that some of the entries of $\frak P^\uparrow_\pi$ might remain
undefined if one of the conditions of type \tht{18} is
not satisfied. Part of the theorem below is that this never happens.

\proclaim{Theorem 10} In the above assumptions, the matrix $\frak
P^\uparrow_\pi$ is well-defined and it is stochastic. Moreover,
$$
\S(\rho_0^+,\dots,\rho_{N-1}^+;\rho_1^-,\dots,\rho_N^-)\, \frak
P^\uparrow_\pi=\S(\rho_0^+,\dots,\rho_{N-1}^+;\rho_1^-,\dots,{}^\sharp\rho_N^-),
$$
where ${}^\sharp{\rho}_N^-=(\rho_N^-,\pi)$. In other words, $\frak
P^\uparrow_\pi$ changes the last specialization of the Schur process
by adding $\pi$ to it.
\endproclaim

The proof of Theorem 10 will be given in Section 9.

Matrices $\frak P^\uparrow_\pi$ describe a certain growth process.
In a similar fashion, one obtains a process of decay. Let us
describe it.

Let $\sigma$ be a nonnegative specialization of $\La$ that `divides'
$\rho_0^+$, that is, there exists a nonnegative specialization
${}^\flat\rho_0^+$ such that $\rho_0^+=({}^\flat\rho_0^+,\sigma)$.
For example, $\sigma$ may coincide with $\rho_0^+$; in that case
${}^\flat\rho_0^+$ is trivial.

Let $\sss$ be the set of pairs of sequences $(\la,\mu)$ as in
\tht{8} with
$$
s_{\la^{(1)}}({}^\flat\rho_0^+)\,s_{\la^{(1)}/\mu^{(1)}}(\rho_1^-)
s_{\la^{(2)}/\mu^{(1)}}(\rho_1^+)\,\cdots
s_{\la^{(N)}/\mu^{(N-1)}}(\rho_{N-1}^+)>0.
$$

Note that if $\sigma=\rho_0^+$ then $\la^{(1)}$ and $\mu^{(1)}$ must
be empty in order for $(\la,\mu)$ to lie in $\sss$.

Define a matrix $\frak P^\downarrow_\sigma$ with rows parameterized
by $\ss$ and columns parameterized by $\sss$ via
$$
\multline \frak
P^\downarrow_\sigma((\la,\mu),(\tilde\la,\tilde\mu))=
P_{\rho_0^+,\sigma}\bigl(\varnothing,\la^{(1)}\uparrow\downarrow
\tilde\la^{(1)}\bigr)\\ \times \prod_{j=1}^{N-1}
P_{\rho_j^-,\sigma}\bigl(\tilde\la^{(j)},\mu^{(j)}\downarrow\downarrow
\tilde\mu^{(j)}\bigr)P_{\rho_j^+,\sigma}\bigl(\tilde\mu^{(j)},\la^{(j+1)}\uparrow\downarrow
\tilde\la^{(j+1)}\bigr).
\endmultline
$$
Notice that the only difference of this definition and that of
$\frak P^\uparrow_\pi$ above, is switching $\pi$ and $\sigma$ and
changing the second arrows from $\uparrow$ to $\downarrow$.

\proclaim{Theorem 11} In the above assumptions, the matrix $\frak
P^\downarrow_\sigma$ is well-defined and it is stochastic. Moreover,
$$
\S(\rho_0^+,\dots,\rho_{N-1}^+;\rho_1^-,\dots,\rho_N^-)\, \frak
P^\downarrow_\sigma=\S({}^\flat\rho_0^+,\dots,\rho_{N-1}^+;\rho_1^-,\dots,\rho_N^-),
$$
where $({}^\flat{\rho}_0^+,\sigma)=\rho_0^+$. In other words, $\frak
P^\downarrow_\sigma$ changes the first specialization of the Schur
process by removing $\sigma$ from it.
\endproclaim

The proof of Theorem 11 will also be given in Section 9.

\example{Remark 12} Both Theorems 10 and 11 can be generalized as
follows. Assume we have an arbitrary sequence of Markov steps of
types $\frak P^\uparrow$ and $\frak P^\downarrow$ applied to an
initial Schur process, and let us denote by $(\lambda(t), \mu(t))$
the result of the application of $t$ first members of the sequence.
One can show that any finite sequence of random partitions of the
form
$$
\multline (\lambda^{(1)}(t_{1,1}),\lambda^{(1)}(t_{1,2}),\dots,
\mu^{(1)}(t_{1,1}'),\mu^{(1)}(t_{1,2}'),\dots,\dots\\
\dots,\mu^{(N-1)}(t_{N-1,1}'),\mu^{(N-1)}(t_{N-1,2}'),\dots,
\lambda^{(N)}(t_{N,1}),\lambda^{(N)}(t_{N,2}),\dots)
\endmultline
$$
forms a Schur process with an explicitly known specializations as
long as
$$
t_{1,1}\ge t_{1,2}\ge \dots\ge t_{1,1}'\ge t_{1,2}'\ge\dots\ge
t_{N-1,1}'\ge t_{N-1,2}'\ge\dots\ge t_{N,1}\ge t_{N,2}\ge\dots,
$$
cf. the last sentence of Section 8.
\endexample

\head 6. Markov chains on the two-sided Schur processes
\endhead

We now aim at formulating (and later proving) a statement for the
two-sided Schur processes that is analogous to Theorem 10.

For two admissible matrices $M_1$ and $M_2$ (`admissible' is
explained in Definition 2), and two signatures $\la$ and $\mu$ of
length $n\ge 1$, we define a probability distribution on
$\Bbb{GT}_n$ (=the set of all signatures of length $n$) via
$$
P_{M_1,M_2}(\la,\mu\,\Vert\,\nu)=const\cdot s_{\nu/\la}(M_1)
s_{\nu/\mu}(M_2),\qquad \nu\in\Bbb{GT}_n.
$$

For an admissible matrix $M$ and a positive number $a$ in the
annulus of analyticity of $H(M;u)$, and for two signatures
$\la\in\Bbb{GT}_{n-1}$ and $\mu\in \Bbb{GT}_n$, we define a
probability distribution on $\Bbb{GT}_n$ via
$$
P_{a,M}(\la,\mu\,\Vert\,\nu)=const\cdot s_{\nu/\la}(a)
s_{\nu/\mu}(M),\qquad \nu\in\Bbb{GT}_n.
$$

In both definitions, we suppose that the set of $\nu$'s giving
nonzero contributions to the right-hand sides is nonempty. Then our
assumptions imply the existence of the normalizing constants.

Similarly to the one-sided Schur process, if $M_1$ and $M_2$ are
both single-$\alpha^\pm$ or single-$\beta^\pm$ matrices, then
$P_{M_1,M_2}$ splits into a product of geometric/Bernoulli random
variables, cf. \tht{5}-\tht{6} and Example 8. For $P_{a,M}$ the same
holds if $M$ is a single-$\alpha^\pm$ or single-$\beta^\pm$ matrix.

Consider the two-sided Schur process of Definition 4, and let
$$
\multline \ss=\Biggl\{(\vec\la^{(1)},\dots,\vec\la^{(N)})\in
(\Bbb{GT}_1)^{c(1)+1}\times\dots\times (\Bbb{GT}_{N})^{c(N)+1} \,\bigl|\\
\prod_{k=1}^N \Bigl( s_{\lambda^{(k,0)}/\lambda^{(k-1,c(k-1))}}(a_k)
\prod_{l=1}^{c(k)}
s_{\lambda^{(k,l)}/\lambda^{(k,l-1)}}\bigl(M^{(k,l)}\bigr)>0\Bigr)\Biggr\},
\endmultline
$$
where $\la^{(0,c(0))}=\varnothing$, cf. \tht{14}. Clearly,
$\operatorname{supp}\T(a_1,\dots,a_N;\M;\Psi)\subset\ss$.

Let $Q$ be an additional admissible matrix such that all the
parameters $a_j$ lie in the analyticity annulus of $H(Q;u)$. Define
a matrix $\frak P_Q$ with rows and columns parameterized by $\ss$
via
$$
\multline
\frak
P_Q((\vec\la^{(1)},\dots,\vec\la^{(N)}),(\vec{\mu}^{(1)},\dots,\vec{\mu}^{(N)}))=\\
 \prod_{k=1}^N
\left(P_{a_k,Q}\left(\lambda^{(k-1,c(k-1))},\lambda^{(k,0)}\,\Vert
\, \mu^{(k,0)}\right) \prod_{l=1}^{c(k)} P_{M^{(k,l)},Q}
\left(\lambda^{(k,l-1)},\lambda^{(k,l)}\,\Vert\,
\mu^{(k,l)}\right)\right)
\endmultline
$$

 The structure of $\frak P_Q$ is such that to compute its
 row indexed by
$(\vec\la^{(1)},\dots,\vec\la^{(N)})$, one first finds $\mu^{(1,0)}$
using $\la^{(1,0)}$, then $\mu^{(1,1)}$ using $\la^{(1,1)}$ and
$\mu^{(1,0)}$, then $\mu^{(1,2)}$ using $\la^{(1,2)}$ and
$\mu^{(1,1)}$, and so on.

\proclaim{Theorem 13} In the above assumptions, the matrix $\frak
P_Q$ is well-defined and it is stochastic. Moreover,
$$
\T(a_1,\dots,a_N;\M;\Psi)\, \frak P_Q=\T(a_1,\dots,a_N;\M;Q\Psi).
$$
\endproclaim

The proof of Theorem 13 will be given in Section 10.

\example{Remark 14} Similarly to Remark 12, a more general statement
can be proved. Assume we have an arbitrary sequence of matrices
${\frak P}_Q$ applied to a two-sided Schur process
$\T(a_1,\dots,a_N;\M;\Psi)$. Denote by
$(\vec\la^{(1)}(t),\dots\vec\la^{(N)}(t))$ the random sequence
obtained after the application of $t$ first matrices. Then any
sequence $\{\la^{(k,l)}(t_{k,l})\}$ forms (a marginal of) an
explicitly describable two-sided Schur process as long as
$(k_1,l_1)\le (k_2,l_2)$ lexicographically implies $t_{k_1,l_1}\ge
t_{k_2,l_2}$.
\endexample

\example{Remark 15} The matrices $\frak P_Q$ are similar to the
growth process defined by $\frak P^\uparrow_\pi$ of the previous
section. One could also define a `decay process' for the two-sided
Schur processes that would be similar to $\frak
P^\downarrow_\sigma$; the application of the corresponding matrix to
$\T(a_1,\dots,a_N;\M;\Psi)$ would reduce $N$ by 1 and remove $a_1$
and $\bigl\{M^{(1,l)}\bigr\}_{l=1}^{c(1)}$ from the set of
parameters.
\endexample

\example{Remark 16} In the setting of Remark 5, one easily shows
that $\frak P^\uparrow_\pi$ and $\frak P_Q$ coincide if
$H(\pi;u)=H(Q;u)$.
\endexample

\head 7. Exact sampling algorithms
\endhead

Let us start with (one-sided) Schur processes. Theorem 10 yields an
exact sampling algorithm that is inductive in $N$.

As the base one can take the empty sequence and $N=0$. Let us
explain the induction step. Assume we already know how to sample
from the Schur process $\Cal
P_{n-1}=\S(\rho_0^+,\dots,\rho_{N-2}^+;\rho_1^-,\dots,\rho_{N-1}^-)$.

Consider the process $\tilde\Cal
P_n=\S(\rho_0^+,\dots,\rho_{N-2}^+,\rho_{N-1}^+;\rho_1^-,\dots,\rho_{N-1}^-,\varnothing)$,
where $\varnothing$ is the trivial specialization. The definition of
the Schur process implies that for this process
$\la^{(N)}=\mu^{(N-1)}=\varnothing$ with probability 1, and the
distribution of the remaining partitions $(\la^{(1)}, \mu^{(1)},\dots,\mu^{(N-2)},
\la^{(N-1)})$ is the same as for $\Cal
P_{n-1}$ that we already know how to sample from by the induction
hypothesis.

In order to obtain a sample of $\Cal
P_n=\S(\rho_0^+,\dots,\rho_{N-2}^+,\rho_{N-1}^+;\rho_1^-,\dots,\rho_{N-1}^-,\rho_N^-)$
we apply the stochastic matrix $\frak P^\uparrow_{\pi}$ with
$\pi=\rho_N^-$ to $\tilde\Cal P_{n}$, cf. Theorem 10. The
application of this matrix requires sequential update from
$\lambda^{(1)}$ up, cf. \tht{19}.

We thus see that if each of
$(\rho_0^+,\dots,\rho_{N-2}^+,\rho_{N-1}^+;\rho_1^-,\dots,\rho_{N-1}^-,\rho_N^-)$
is a single-$\alpha$ or a single-$\beta$ specializations (or
trivial), then exact sampling is reduced to sampling a finite number
of independent geometric/Bernoulli random variables. Noting that in
the algorithm for the $N$th step one does not have to use a single
$\frak P^\uparrow_{\pi}$ with $\pi=\rho_N^-$, but can instead use a
sequence of $\frak P^\uparrow_{\pi_k}$ with
$\rho_N^-=(\pi_1,\pi_2,\dots)$, we see that the a similar
reduction holds for the Schur processes with all specializations having finitely
many nonzero $\alpha$'s and $\beta$'s (and $\gamma=0$).

For the measures $q^{volume}$ on skew plane partitions considered in
Section 3, the algorithm can be implemented as follows (we use
Section 3 and Example 8 below).
\medskip
{\tt Initiate by assigning
$\la^{(1)}=\dots=\la^{(A+B)}=\varnothing$.

For $k$ running from 2 to $(A+B)$

    \quad If $k\notin\Cal L(\pi)$ then

      \quad \quad For l running from 1 to $(k-1)$

         \quad\quad\quad If $l\in\Cal L(\pi)$ then
         $\lambda^{(l+1)}:=\nu$ with
         $\nu$ distributed as

         \quad\quad\quad\quad $G_{\max\{\la_1^{(l)},\la_1^{(l+1)}\},+\infty}^{q^{k-l}}(\nu_1)
                                \prod\limits_{j\ge 2}
                                G_{\max\{\la_j^{(l)},\la_j^{(l+1)}\},\min\{\la_{j-1}^{(l)},\la_{j-1}^{(l+1)}\}}^{q^{k-l}}(\nu_j)$

         \quad\quad\quad If $l\notin\Cal L(\pi)$ then
         $\lambda^{(l+1)}:=\nu$ with
         $\nu$ distributed as

         \quad\quad\quad\quad $G^{q^{k-l}}_{\max\{\la_{2}^{(l)},\la_1^{(l+1)}\},\la_1^{(l+1)}}(\nu_1)
         \prod\limits_{j\ge 2}
G^{q^{k-l}}_{\max\{\la_{j+1}^{(l)},\la_{j}^{(l+1)}\},\min\{\la_{j}^{(l)},\la_{j-1}^{(l+1)}\}}(\nu_j)$

\quad \quad End of $l$-cycle

End of $k$-cycle }
\smallskip

At the end of each $k$-step we see an exact sample of the measure
$q^{volume}$ on plane partitions with a smaller support. The number
of nontrivial one-dimensional samples needed to go through the
$k$-step with $k\notin\Cal L(\pi)$ is the number of boxes in this
support. It is not difficult to see that this number is at most $A$
for the smallest $k\notin\Cal L(\pi)$, it is at most $2A$ for the
next one and so on, so that the total number of one-dimensional
samples needed is at most $AB(B+1)/2$. The maximum is achieved at
$\Cal L(\pi)=\{1,\dots,A\}$, i.e. when the plane partitions are
supported by the full $A\times B$ box.

%\vskip -0.0cm\hskip -0.5cm
% {\pdfximage width 60 mm
%{cusp2.jpg}\pdfrefximage \pdflastximage}

%\vskip -5.0cm\hskip 5.5cm {\pdfximage width 60 mm
%{cusp_avg1.jpg}\pdfrefximage \pdflastximage}

$$ \epsffile{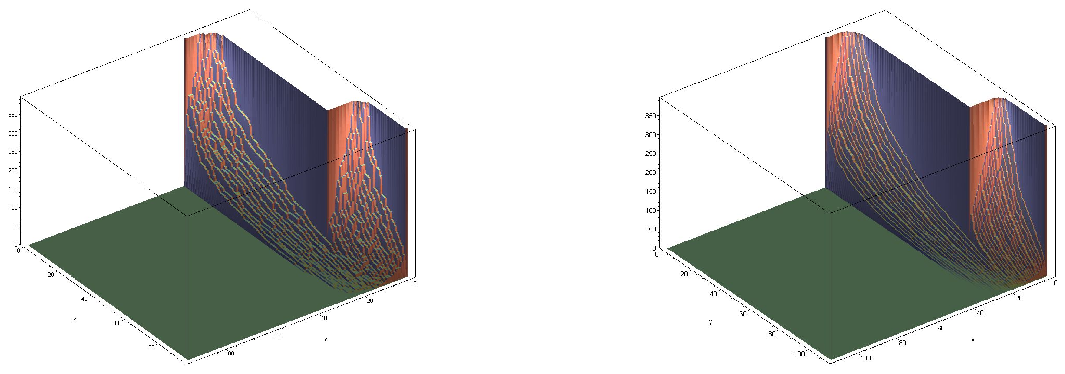} $$

The above figures show a sample for a specific back wall profile,
and an average over ten samples with the same back wall. A limit
shape and its cusp are clearly visible, cf. \cite{OR2}.

Finally, note that a very similar algorithm would sample skew plane
partitions with weights of the form $\prod q_j^{|\la^{(j)}|}$.

Let us now discuss the two-sided Schur process. First, let us
restrict ourselves to the case when $\Psi$ is Toeplitz. Then if all
$H(M^{(k,l)};u^{-1})$ and $\psi(u)$ are analytic in a disc of radius $>1$
(not just in an annulus containing the unit circle), then the
two-sided Schur process lives on signatures with nonnegative
coordinates and it constitutes a special case of the (one-sided)
Schur process, cf. Remark 5. Consequently, if all $M^{(k,l)}$ and
$\Psi$ are admissible matrices with $M^{(k,l)}$ having finitely many
$\alpha^-$ and $\beta^-$ nonzero parameters (all others are zero),
and $\Psi$ having finitely many $\alpha^+$ and $\beta^+$ nonzero
parameters, the inductive algorithm for the Schur process described
above reduces sampling to a finite number of independent samples of
geometric/Bernoulli random variables.

On the other hand, Theorem 13 allows us to add finitely many
$\alpha^\pm$ and $\beta^\pm$ parameters to $\Psi$ by sampling from
independent geometric/Bernoulli distributions. Hence, we can relax
the assumption on $\Psi$ in the previous paragraph by requiring that
it has finitely many $\alpha^\pm$ and $\beta^\pm$ parameters.

%\vskip -0.0cm\hskip -0.5cm
% {\pdfximage width 65 mm
%{gt_large1.jpg}\pdfrefximage \pdflastximage}

%\vskip -5.9cm\hskip 5.5cm {\pdfximage width 65 mm
%{gt_large1_avg.jpg}\pdfrefximage \pdflastximage}

$$ \epsffile{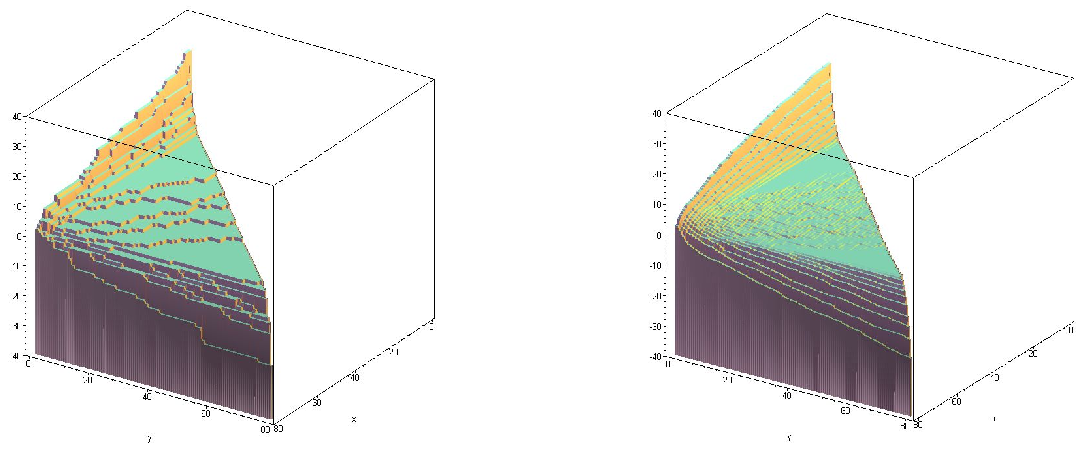} $$

The figure above shows a sample of the path measure and the average
over ten samples for the extreme character of $U(\infty)$ with
$$
\alpha_1^+=\dots\alpha_{10}^+=\frac 1{10},\quad
\beta_1^+=\dots\beta_5^+=\frac12,\quad\alpha_1^-=\dots=\alpha_{10}^-=\frac
1{10},
$$
and all other parameters being zero, cf. Section 4. The first order
asymptotic behavior of such measures as the path length goes to
infinity and parameters remain fixed is known, see \cite{OO}.

\head 8. A general construction of multivariate Markov chains
\endhead

The general construction of this section will be used in Sections 9
and 10 to prove Theorems 10, 11, and 13.

Let $(\SS_1,\dots,\SS_n)$ and $(\tilde\SS_1,\dots,\tilde\SS_n)$ be
two $n$-tuples of discrete countable sets, $P_1,\dots,P_n$ be
stochastic matrices defining Markov chains $\SS_j\to\tilde\SS_j$.
Also let $\Lambda_1^2,\dots,$ $\Lambda_{n-1}^n$ and
$\tilde\Lambda_1^2,\dots,\tilde\Lambda_{n-1}^n$ be stochastic links
between these sets:
$$
\aligned P_k:\SS_k\times\tilde\SS_k\to [0,1],\quad
\sum_{y\in\tilde\SS_k}P_k(x,y)=1,\quad x\in \SS_k,\quad k=1,\dots,n;\\
\Lambda_{k-1}^k:\SS_k\times\SS_{k-1}\to [0,1],\quad
\sum_{y\in\SS_{k-1}}\Lambda_{k-1}^k(x,y)=1,\quad x\in \SS_k,\quad
k=2,\dots,n;
\\
\tilde\Lambda_{k-1}^k:\tilde\SS_k\times\tilde\SS_{k-1}\to
[0,1],\quad
\sum_{y\in\tilde\SS_{k-1}}\tilde\Lambda_{k-1}^k(x,y)=1,\quad x\in
\tilde\SS_k,\quad k=2,\dots,n.
\endaligned
$$

Assume that these matrices satisfy the commutation relations
$$
\Delta^k_{k-1}:=\Lambda^k_{k-1}P_{k-1}=P_k\tilde\Lambda^k_{k-1},\qquad
k=2,\dots,n. \tag 20
$$

We will define a {\it multivariate\/} Markov chain $P^{(n)}$ between
the state spaces
$$
\SS^{(n)}=\Bigl\{(x_1,\dots,x_n)\in\SS_1\times\cdots\times\SS_n\mid
\prod_{k=2}^n\Lambda_{k-1}^k(x_k,x_{k-1})\ne 0\Bigr\}
$$
and
$$
\tilde\SS^{(n)}=\Bigl\{(x_1,\dots,x_n)\in\tilde\SS_1\times\cdots\times\tilde\SS_n\mid
\prod_{k=2}^n\tilde\Lambda_{k-1}^k(x_k,x_{k-1})\ne 0\Bigr\}.
$$
 The transition probabilities for the Markov chain
$P^{(n)}$ are defined as (we use the notation $X_n=(x_1,\dots,x_n)$,
$Y_n=(y_1,\dots,y_n)$)
$$
P^{(n)}(X_n,Y_n)=
P_1(x_1,y_1)\prod\limits_{k=2}^n\dfrac{P_k(x_k,y_k)\tilde\Lambda_{k-1}^k(y_k,y_{k-1})}
{\Delta^k_{k-1}(x_k,y_{k-1})} \tag 21
$$
if $\prod_{k=2}^n\Delta^k_{k-1}(x_k,y_{k-1})>0$, and $0$ otherwise.

One way to think of $P^{(n)}$ is as follows.

Starting from $X=(x_1,\dots,x_n)$, we first choose $y_1$ according
to the transition matrix $P_1(x_1,y_1)$, then choose $y_2$ using
$\frac{P_2(x_2,y_2)\tilde\Lambda_{1}^2(y_2,y_1)}{\Delta^2_{1}(x_2,y_1)}$,
which is the conditional distribution of the middle point in the
successive application of $P_2$ and $\tilde\Lambda^2_1$ provided
that we start at $x_2$ and finish at $y_1$, after that we choose
$y_3$ using the conditional distribution of the middle point in the
successive application of $P_3$ and $\tilde\Lambda^3_2$ provided
that we start at $x_3$ and finish at $y_2$, and so on. Thus, one
could say that $Y$ is obtained from $X$ by the sequential update.

\proclaim{Proposition 17} Let $m_n$ be a probability measure on
$\SS_n$. Let $m^{(n)}$ be a probability measure on $\SS^{(n)}$
defined by
$$
m^{(n)}(X_n)=m_n(x_n)\Lambda^n_{n-1}(x_n,x_{n-1})\cdots
\Lambda^2_1(x_2,x_1),\qquad X_n=(x_1,\dots,x_n)\in\SS^{(n)}.
$$
Set $\tilde m_n=m_nP_n$ and
$$
\tilde m^{(n)}(X_n)=\tilde
m_n(x_n)\tilde\Lambda^n_{n-1}(x_n,x_{n-1})\cdots
\tilde\Lambda^2_1(x_2,x_1),\qquad
X_n=(x_1,\dots,x_n)\in\tilde\SS^{(n)}.
$$
Then $m^{(n)}P^{(n)}=\tilde m^{(n)}$.
\endproclaim
\demo{Proof} The argument is straightforward. Indeed,
$$
\multline
m^{(n)}P^{(n)}(Y_n)=\sum_{X_n\in\SS^{(n)}}m_n(x_n)\Lambda^n_{n-1}(x_n,x_{n-1})\cdots
\Lambda^2_1(x_2,x_1)\\
\times
P_1(x_1,y_1)\prod\limits_{k=2}^n\dfrac{P_k(x_k,y_k)\tilde\Lambda_{k-1}^k(y_k,y_{k-1})}
{\Delta^k_{k-1}(x_k,y_{k-1})}\,.
\endmultline
$$
Extending the sum to $x_1\in\SS_1$ adds 0 to the right-hand side.
Then we can use  relation \tht{20} to compute the sum over $x_1$,
removing $\Lambda^2_1(x_2,x_1)$, $P_1(x_1,y_1)$ and
$\Delta^2_1(x_2,y_1)$ from the expression. Similarly, we sum
consecutively over $x_2,\dots,x_n$, and this gives the needed
result.\qed
\enddemo

Proposition 17 will be used to prove Theorems 10, 11, and 13. A more
general \cite{BF, Proposition 2.7} is needed to prove the statements
mentioned in Remarks 12 and 14.

\head 9. Application to the Schur processes
\endhead

In this section we prove Theorems 10 and 11.

Let us start by putting the Schur process
$\S(\rho_0^+,\dots,\rho_{N-1}^+;\rho_1^-,\dots,\rho_N^-)$ of
Definition 3 into the framework of the previous section.

We need some general definitions.

Let $y,z,t$ be nonnegative specializations of $\Lambda$. Set
$$
\align p_{\lambda\mu}^{\uparrow}(y;z)&=\frac 1{H(y;z)}
\frac{s_{\mu}(y)}{s_{\la}(y)}\,s_{\mu/\la}(z),\qquad \la,\mu\in \Y(y),\\
p_{\lambda\nu}^{\downarrow}(y;t)&=
\frac{s_{\nu}(y)}{s_{\la}(y,t)}\,s_{\la/\nu}(t),\qquad\qquad
\la\in\Y(y,t),\ \nu\in\Y(y),
\endalign
$$
where $\Y(\rho)=\{\kappa\in\Y\mid s_\kappa(\rho)>0\}$, and for the
first definition we assume that $H(y;z)=\sum_{\kappa\in\Y}
s_\kappa(y)s_\kappa(z)<\infty$.

Relations \tht{11} and \tht{12} imply that the matrices
$$
p^\uparrow(y;z)=\bigl[p_{\lambda\mu}^{\uparrow}(y;z)\bigr]_{\la,\mu\in\Y(y)}\quad
\text{and}\quad p^\downarrow(y;t)=
\bigl[p_{\lambda\nu}^{\downarrow}(y;t)\bigr]_{\la\in\Y(y,t),\nu\in\Y(y)}
$$
are stochastic:
$$
\sum_{\mu\in\Y(y)} p_{\lambda\mu}^{\uparrow}(y;z)=\sum_{\nu\in\Y(y)}
p_{\lambda\nu}^{\downarrow}(y;t)=1.
$$

It is immediate to see that $p^\uparrow$ and $p^\downarrow$ act well
on the Schur measures:
$$
\S(x;y)p^{\uparrow}(y;z)=\S(x,z;y),\qquad
\S(x;y,t)p^{\downarrow}(y;t)=\S(x;y). \tag 22
$$
Observe that $\S(\rho_1;\rho_2)=\S(\rho_2;\rho_1)$, so the
parameters of the Schur measures in these relations can also be
permuted.

\proclaim{Proposition 18} Let $y,z,z_1,z_2,t_1,t_2$ be nonnegative
specializations of $\La$. Then we have the commutativity relations
$$
\align p^{\uparrow}(y;z_1)p^{\uparrow}(y;z_2)&=
p^{\uparrow}(y;z_2)p^{\uparrow}(y;z_1),\\
p^{\downarrow}(y,t_2;t_1)p^{\downarrow}(y;t_2)&=
p^{\downarrow}(y,t_1;t_2)p^{\downarrow}(y;t_1),\\
p^{\uparrow}(y,t;z)p^{\downarrow}(y;t)&=
p^{\downarrow}(y;t)p^{\uparrow}(y;z),
\endalign
$$
where for the first relation we assume $H(y;z_1,z_2)<\infty$, and
for the third relation we assume $H(y,t;z)<\infty$.
\endproclaim

\demo{Proof} The arguments for all three identities are similar; we only
give the proof of the third one which is in a way the
hardest. We have
$$
\multline \sum_{\mu}
p_{\la\mu}^{\uparrow}(y,t;z)p^{\downarrow}_{\mu\nu}(y;t)=
\frac{1}{H(y,t;z)}\sum_{\mu\in\Y(y,t)}
\frac{s_\mu(y,t)}{s_\la(y,t)}s_{\mu/\la}(z)\,\frac{s_\nu(y)}{s_\mu(y,t)}s_{\mu/\nu}(t)\\=
\frac{1}{H(y,t;z)}\frac{s_\nu(y)}{s_\la(y,t)}\sum_{\mu\in\Y}
s_{\mu/\la}(z)s_{\mu/\nu}(t)=\frac{H(t;z)}{H(y,t;z)}
\frac{s_\nu(y)}{s_\la(y,t)} \sum_{\kappa\in\Y}
s_{\la/\kappa}(t)s_{\nu/\kappa}(z)\\=
\frac{1}{H(y;z)}\sum_{\kappa\in\Y(y)}
\frac{s_\kappa(y)}{s_\la(y,t)}s_{\la/\kappa}(t)\,\frac{s_\nu(y)}{s_\kappa(y)}s_{\nu/\kappa}(z)=
\sum_{\kappa\in\Y(y)}
p_{\la\kappa}^{\downarrow}(y;t)p^{\uparrow}_{\kappa\nu}(y;z),
\endmultline
$$
where along the way we extended the summation in $\mu$ from
$\Y(y,t)$ to $\Y$ because $s_\nu(y)s_{\mu/\nu}(t)>0$ implies
$s_\mu(y,t)>0$ by \tht{12}; we used \tht{11} to switch from $\mu$ to
$\kappa$, and finally we restricted the summation in $\kappa$ from
$\Y$ to $\Y(y)$ because $s_{\nu}(y)s_{\nu/\kappa}(z)>0$ implies
$\kappa\subset\nu$ and $s_{\kappa}(y)>0$.\qed
\enddemo

We are now ready to return to the Schur process
$\S(\rho_0^+,\dots,\rho_{N-1}^+;\rho_1^-,\dots,\rho_N^-)$.

Set $n=2N-1$ and
$$
\align &\SS_{2j-1}=\Y\bigl(\rho^+_{[0,j-1]}\bigr), \qquad
j=1,\dots,N;\\
&\SS_{2k}=\Y\bigl(\rho^+_{[0,k-1]}\bigr),\qquad k=1,\dots,N-1.
\endalign
$$
Since $\la^{(j)}$ and $\mu^{(k)}$ are distributed according to the
Schur measures $\S(\rho^+_{[0,j-1]};\rho^-_{[j,N]})$ and
$\S(\rho^+_{[0,k-1]};\rho^-_{[k+1,N]})$ respectively, the
projections of the support of the Schur process to these coordinates
lie inside $\SS_{2j-1}$ and $\SS_{2k}$, respectively.

Define the stochastic links by
$$
\align
&\Lambda^{2j+1}_{2j}=p^{\downarrow}(\rho^+_{[0,j-1]};\rho_{j}^+),\qquad
j=1,\dots,N-1;\\
&\Lambda^{2j}_{2j-1}=p^\uparrow(\rho^+_{[0,j-1]};\rho_j^-),\qquad
j=1,\dots,N-1.
\endalign
$$

 One immediately verifies the formula
$$
\multline
\S(\rho_0^+,\dots,\rho_{N-1}^+;\rho_1^-,\dots,\rho_N^-)(\la,\mu)\\=
\S\bigl(\rho_{[0,N-1]}^+;\rho_{N}^-\bigr)(\la^{(N)})\prod_{k=1}^{N-1}
\biggl(\Lambda^{2k+1}_{2k}\bigl(\lambda^{(k+1)},\mu^{(k)}\bigr)
\,\Lambda^{2k}_{2k-1}\bigl(\mu^{(k)},\la^{(k)}\bigr)\biggr),
\endmultline
\tag 23
$$
cf. the definition of $m^{(n)}$ in Proposition 17.

\demo{Proof of Theorem 10} We apply Proposition 17. Set
$\tilde\SS_j=\SS_j$ for $j=1,\dots,n$,
$\tilde\La^{j}_{j-1}=\La^{j}_{j-1}$ for $j=2,\dots,n$, and also
$$
\align m_n&=\S(\rho_{[0,N-1]}^+;\rho_{N}^-),\\
P_{2j-1}&=p^{\uparrow}(\rho^+_{[0,j-1]};\pi),\quad
j=1,\dots,N;\\
P_{2j}&=p^{\uparrow}(\rho^+_{[0,j-1]};\pi),\quad j=1,\dots,N-1.
\endalign
$$
The commutation relations \tht{20} follow from Proposition 18, and
the matrix of transition probabilities $P^{(n)}$ from \tht{21} is
easily seen to coincide with $\frak P^\uparrow_\pi$. The claim now
follows from \tht{23}, Proposition 17, and the relation (cf.
\tht{22})
$$
\S\bigl(\rho_{[0,N-1]}^+;\rho_{N}^-\bigr)P_n=\S(\rho_{[0,N-1]}^+;\rho_{N}^-,\pi).\qed
$$
\enddemo

\demo{Proof of Theorem 11} We also apply Proposition 17. This time
we need to modify the state spaces:
$$
\align
&\tilde\SS_{2j-1}=\Y\bigl({}^\flat\rho^+_0,\rho^+_{[1,j-1]}\bigr),
\qquad
j=1,\dots,N;\\
&\tilde\SS_{2k}=\Y\bigl({}^\flat\rho^+_0,\rho^+_{[1,k-1]}\bigr),\qquad
k=1,\dots,N-1.
\endalign
$$
Also set
$$
\align
&\tilde\Lambda^{2j+1}_{2j}=p^{\downarrow}({}^\flat\rho_0^+,\rho^+_{[1,j-1]};\rho_{j}^+),\qquad
j=1,\dots,N-1;\\
&\tilde\Lambda^{2j}_{2j-1}=p^\uparrow({}^\flat\rho_0^+,\rho^+_{[1,j-1]};\rho_j^-),\qquad
j=1,\dots,N-1;
\endalign
$$
 and
$$
\align m_n&=\S(\rho_{[0,N-1]}^+;\rho_{N}^-),\\
P_{2j-1}&=p^{\downarrow}({}^\flat\rho_0^+,\rho^+_{[1,j-1]};\sigma),\qquad
j=1,\dots,N;\\
P_{2j}&=p^{\downarrow}({}^\flat\rho_0^+,\rho^+_{[1,j-1]};\sigma),\qquad
j=1,\dots,N-1,
\endalign
$$
Again, the commutation relations \tht{20} follow from Proposition
18, and the matrix of transition probabilities $P^{(n)}$ from
\tht{21} coincides with $\frak P^\downarrow_\sigma$. The claim
follows from \tht{23}, Proposition 17, and the relation (cf.
\tht{22})
$$
\S\bigl(\rho_{[0,N-1]}^+;\rho_{N}^-\bigr)P_n=\S({}^\flat\rho_0^+,\rho_{[1,N-1]}^+;\rho_{N}^-).\qed
$$
\enddemo

\head 10. Application to the two-sided Schur processes\endhead

Let us put the two-sided Schur process of Definition 4 into the
general framework.

We need some notation. For $n\ge 1$, an admissible matrix $M$, cf.
Definition 2, and $a_1,\dots,a_n> 0$ in the analyticity annulus of
$H(M;u^{-1})$, define
$$
T_{\la\mu}(a_1,\dots,a_n;M)=\frac 1{\prod\limits_{j=1}^n
H(M;a_j^{-1})}
\frac{\det\bigl[a_i^{\mu_j-j}\bigr]_{i,j=1}^n}{\det\bigl[a_i^{\la_j-j}\bigr]_{i,j=1}^n}\,s_{\la/\mu}(M),
\qquad \la,\mu\in\Bbb{GT}_n.
$$
For arbitrary $a_1,\dots,a_n> 0$ also set ($\la\in\Bbb{GT}_n$,
$\mu\in\Bbb{GT}_{n-1}$)
$$
T_{\la\mu}(a_1,\dots,a_n)=\frac 1{a_n}{\prod\limits_{j=1}^{n-1}
\left( \frac1{a_n}-\frac{1}{a_j}\right)}
\frac{\det\bigl[a_i^{\mu_j-j}\bigr]_{i,j=1}^{n-1}}{\det\bigl[a_i^{\la_j-j}\bigr]_{i,j=1}^n}\,s_{\la/\mu}(a_n).
$$
Thus, we have matrices $T(a_1,\dots,a_n;M)$ with rows and column
parameterized by $\Bbb{GT}_n$, and matrices $T(a_1,\dots,a_n)$ with
rows parameterized by $\Bbb{GT}_{n}$ and columns parameterized by
$\Bbb{GT}_{n-1}$.

\proclaim{Proposition 19} In the above assumptions, the matrices
$T(a_1,\dots,a_n;M)$ and $T(a_1,\dots,a_n)$ are stochastic, and the
following commutation relation holds:
$$
T(a_1,\dots,a_n;M)\,T(a_1,\dots,a_n)=T(a_1,\dots,a_n)\,T(a_1,\dots,a_{n-1};M).
$$
For admissible matrices $M_1,M_2$ and $a_1,\dots,a_n>0$ in the
analyticity annuli of $H(M_i;u^{-1})$, $i=1,2$, we also have the
commutation relation
$$
T(a_1,\dots,a_n;M_1)\,T(a_1,\dots,a_n;M_2)=T(a_1,\dots,a_n;M_2)\,T(a_1,\dots,a_{n};M_1).
$$

\endproclaim
\demo{Proof} Follows from Propositions 2.8-2.10 and Lemma 2.13(ii)
of \cite{BF}. \qed
\enddemo

Consider now the two-sided Schur process $\T(a_1,\dots,a_N;\M;\Psi)$
of Definition 4. Set $n=c(1)+\dots+c(N)+N$, and ($c(0):=0$)
$$
\SS_j=\Bbb{GT}_k,\qquad  c(k-1)+k\le j\le c(k)+k,\quad k=1,\dots,N.
$$
Define the stochastic links by
$$
\align
\Lambda^{c(k-1)+k}_{c(k-1)+k-1}&=T(a_1,\dots,a_{k}),\qquad\qquad\quad
k=2,\dots,N;\\
\Lambda^{c(k-1)+k+l}_{c(k-1)+k+l-1}
&=T(a_1,\dots,a_{k};M^{(k,l)}),\qquad k=1,\dots,N,\quad
l=1,\dots,c(k).
\endalign
$$
Also define a probability distribution $m_n^\Psi$ on
$\SS_n=\Bbb{GT}_{N}$ via
$$
m_n^\Psi(\la)=\frac{\det\bigl[a_i^{\la_j-j}\bigr]_{i,j=1}^N\det\bigl[\Psi_{\la_i-i,-j}\bigr]_{i,j=1}^N}{\det\bigl[a_i^{-j}\Psi_j(a_i)\bigr]_{i,j=1}^N}\,,
\qquad \la\in\Bbb{GT}_N,
$$
where we used the notation \tht{15}.

These definitions imply that
$$
\multline
\T(a_1,\dots,a_N;\M;\Psi)(\vec\la^{(1)},\dots,\vec\la^{(N)})\\=
m_n^\Psi\bigl(\la^{(N,c(N))}\bigr)\La^n_{n-1}
\bigl(\la^{(N,c(N))},\la^{(N,c(N-1))}\bigr)\dots\La^2_1\bigl(\la^{(1,1)},\la^{(1,0)}\bigr).
\endmultline
$$
Note that this proves formula \tht{16} for the partition function
since
$$
\det\left[a_i^{-j}\right]_{i,j=1}^N=\prod_{k=1}^n \frac
1{a_k}\prod_{j=1}^{k-1}\left(\frac 1{a_k}-\frac{1}{a_j}\right).
$$
 \demo{Proof of Theorem 13} Once again we apply Proposition 17. We
set $\tilde\SS_j=\SS_j$ for $j=1,\dots,n$;
$\tilde\La^j_{j-1}=\La^j_{j-1}$ for $j=2,\dots,n$; and
$m_n=m_n^\Psi$,
$$
P_j=T(a_1,\dots,a_k;Q^t),\qquad c(k-1)+k\le j\le c(k)+k,\quad
k=1,\dots,N.
$$
Note that $H(Q^t;u)=H(Q;u^{-1})$ and
$s_{\la/\mu}(Q^t)=s_{\mu/\la}(Q)$ for signatures $\la$ and $\mu$ of
the same length.

The claim now follows from Proposition 17 as the needed
commutativity relations are given in Proposition 19, and by the
Cauchy-Binet identity
$$
\multline (m_n^\Psi P_n)(\mu)=\frac
1{\det\bigl[a_i^{-j}\Psi_j(a_i)\bigr]_{i,j=1}^N}\frac1{\prod_{j=1}^n
H(Q;a_j)}\\ \times\sum_{\la\in\Bbb{GT}_N}
{\det\bigl[a_i^{\la_j-j}\bigr]_{i,j=1}^N\det\bigl[\Psi_{\la_i-i,-j}\bigr]_{i,j=1}^N}
\frac{\det\bigl[a_i^{\mu_j-j}\bigr]_{i,j=1}^n}{\det\bigl[a_i^{\la_j-j}\bigr]_{i,j=1}^n}\,s_{\mu/\la}(Q)
=m_n^{Q\Psi}(\mu).
\endmultline
$$
\enddemo

\enddocument
\end